\documentclass{amsart}
\usepackage{amsmath, amssymb}
\usepackage{amsfonts}
\usepackage[arrow,matrix,curve,cmtip,ps]{xy}

\usepackage{amsthm}

\allowdisplaybreaks

\newtheorem{theorem}{Theorem}[section]
\newtheorem{lemma}[theorem]{Lemma}
\newtheorem{proposition}[theorem]{Proposition}
\newtheorem{corollary}[theorem]{Corollary}

\newtheorem*{theorem*}{Theorem}
\theoremstyle{remark}
\newtheorem{remark}[theorem]{Remark}
\newtheorem{definition}[theorem]{Definition}
\newtheorem{example}[theorem]{Example}

\numberwithin{equation}{section}

\newcommand{\Z}{\mathbb{Z}}
\newcommand{\N}{\mathbb{N}}

\newcommand{\C}{\mathbb{C}}
\newcommand{\T}{\mathbb{T}}
\newcommand{\K}{\mathcal{K}}

\newcommand{\Li}{\mathcal{L}}
\newcommand{\OA}{\mathcal{O}_A}
\newcommand{\OX}{\mathcal{O}_X}
\newcommand{\Gr}{\mathrm{Gr}}

\newcommand{\G}{\mathcal{G}}
\newcommand{\F}{\mathcal{F}}
\newcommand{\Ri}{\mathfrak{R}}
\newcommand{\Pow}{\mathcal{P}}
\newcommand{\Hi}{\mathcal{H}}

\newcommand{\coker}{\operatorname{coker }}

\newcommand{\rank}{\operatorname{rank}}

\begin{document}
\title{Simplicity of ultragraph algebras}

\author{Mark Tomforde %\thanks{support}\\
		%Department of Mathematics\\
		%Dartmouth College\\
		%Hanover, NH  03755\\
		%mark.tomforde@dartmouth.edu
}

\address{Department of Mathematics\\ Dartmouth College\\
Hanover\\ NH 03755-3551\\ USA}

\curraddr{Department of Mathematics\\ University of Iowa\\
Iowa City\\ IA 52242\\ USA}

\email{tomforde@math.uiowa.edu}

%\thanks{This research was supported by }

\date{\today}
\subjclass{46L55}

\keywords{$C^*$-algebras, simple $C^*$-algebras,
ultragraphs, ultragraph algebras, graph algebras, Exel-Laca
algebras, purely infinite, AF-algebras}

\begin{abstract}

In this paper we analyze the structure of $C^*$-algebras
associated to ultragraphs, which are generalizations of
directed graphs.  We characterize the simple ultragraph algebras as well as deduce necessary and sufficient
conditions for an ultragraph algebra to be purely infinite
and to be AF.  Using these techniques we also produce an
example of an ultragraph algebra that is neither a graph
algebra nor an Exel-Laca algebra.  We conclude by proving
that the $C^*$-algebras of ultragraphs with no sinks are
Cuntz-Pimsner algebras.
\end{abstract}

\maketitle

%%%%%%%%%%%%%%%%%%%%%%%%%%%%%%%%%%%%%%%%%%%%%%%%%%%%%%%%%%%%%%%
\section{Introduction}
%%%%%%%%%%%%%%%%%%%%%%%%%%%%%%%%%%%%%%%%%%%%%%%%%%%%%%%%%%%%%%%

In \cite{Tom3} a generalization of a directed graph,
called an ultragraph, was defined.  In analogy with
the $C^*$-algebras of directed graphs, it was also shown how
to associate a $C^*$-algebra $C^*(\G)$ to an ultragraph
$\G$.  These ultragraph algebras include the
$C^*$-algebras of graphs \cite{KPRR,KPR,BPRS,FLR} as well as
the Exel-Laca algebras of \cite{EL}.  Furthermore, it
was shown that many of the techniques used for graph
algebras can be applied to obtain similar results for
ultragraph algebras.  This has many important
consequences.  First, one can now study Exel-Laca algebras
in terms of ultragraphs.  Thus the frequently complicated
and cumbersome matrix manipulations involved in studying
Exel-Laca algebras may be replaced by graphical techniques
that are often easier to deal with as well as more visual. 
In addition, since the classes of graph algebras and
Exel-Laca algebras each contain $C^*$-algebras that are not
in the other, similar results concerning the two classes
have often had to be proven separately for each class. 
Because ultragraph algebras contain both of these
classes, they provide a context in which these similar
results can be proven once and then applied to the special
cases of graph algebras and Exel-Laca algebras.  

In this paper we build upon the work in \cite{Tom3} and
analyze the structure of $C^*(\G)$.  Throughout we have two
goals.  First, we wish to show that graph algebra techniques
can be used to obtain many results concerning $C^*(\G)$ and
that many properties of $C^*(\G)$ can be read off from the
ultragraph $\G$.  Second, we wish to convince the
reader that the ultragraph approach provides a more
convenient method for studying Exel-Laca algebras.  In their
seminal paper \cite{EL}, Exel and Laca describe how
to associate a graph $\Gr (A)$ to a $\{0,1\}$-matrix $A$. 
Throughout their analysis many conditions are
stated in terms of the graph $\Gr(A)$ and it is shown that
certain properties of $\OA$ are reflected in $\Gr (A)$. 
As in \cite{Tom3} we shall associate an ultragraph $\G_A$
to $A$ for which $C^*(\G_A)$ is canonically isomorphic to
$\OA$.  We shall show that the ultragraph $\G_A$
provides much of the same information as $\Gr (A)$, and
in addition there are aspects of $\OA$ that can be easily
obtained from $\G_A$ but not from $\Gr
(A)$.  In particular, we examine how the simplicity of $\OA$
is reflected in $\G_A$.

After some preliminaries, we begin in
\S\ref{simp-sec} by considering the ideals of $C^*(\G)$ and
determining necessary and sufficient conditions for
$C^*(\G)$ to be simple.  Finding conditions for simplicity in
graph algebras and Exel-Laca algebras has been an elusive
goal of many authors in the past few years.  It was not
until recently that such conditions were obtained, and
the preliminary work involved many partial
results as well as high-powered techniques and sophisticated
tools.  Building on the simplicity criteria for
Cuntz-Krieger algebras \cite[Theorem~2.14]{CK} conditions
for simplicity of $C^*$-algebras of certain graphs
were obtained in \cite[Corollary 6.8]{KPRR} and similar
results for row-finite graphs were obtained in
\cite[Proposition 5.1]{BPRS}.  In
\cite{FLR} $C^*$-algebras of arbitrary (i.e., not necessarily
row-finite) graphs were introduced and it was shown that
transitivity of the graph was a sufficient (but not
necessary) condition for simplicity of the $C^*$-algebra 
\cite[Theorem~3]{FLR}.  In
\cite[Corollary~4.5]{FR} it was shown that for graphs in
which every vertex emits infinitely many edges, transitivity
was also a necessary condition for simplicity.  In addition,
Exel and Laca gave sufficient conditions for simplicity
of the Exel-Laca algebras in \cite[Theorem~14.1]{EL}. 
Necessary and sufficient conditions for simplicity of
Exel-Laca algebras were finally obtained by Szyma\'nski in
\cite[Theorem~8]{Szy} and his result
could be adapted to give necessary and
sufficient conditions for simplicity of
$C^*$-algebras of arbitrary graphs \cite[Theorem~12]{Szy}. 
His conditions for the Exel-Laca algebras $\OA$ were stated
in terms of saturated hereditary subsets of the index set
of $A$, and his conditions for graph algebras were stated in
terms of saturated hereditary subsets of the graph's
vertices.  Shortly afterwards independent results of
\cite[Theorem~4]{Pat} and \cite[Corollary~2.14]{DT} also
gave necessary and sufficient conditions for simplicity of
graph algebras in terms of reachability of certain vertices
in the graph.  

In this paper we give necessary and sufficient conditions
for an ultragraph algebra to be simple.  We state this
result in two ways.  In Theorem~\ref{Woj-simplicity} we
give the result in terms of saturated hereditary
subcollections, and as one would expect the result is very
much like that of Szyma\'nski's in  \cite[Theorem~12]{Szy}. 
In addition, in Theorem~\ref{DT-simplicity} we give a
characterization of simplicity in terms of reachability of
certain vertices.  Although this result contains
\cite[Theorem~4]{Pat} and \cite[Corollary~2.14]{DT} as
special cases, it is a much less obvious generalization.  We
conclude \S\ref{simp-sec} with an example showing that the
ultragraph $\G_A$ is a better tool than the graph $\Gr(A)$
for determining the simplicity of the Exel-Laca algebra
$\OA$.

In \S\ref{AF-pi} we give necessary and sufficient conditions
for $C^*(\G)$ to be purely infinite and to be AF. 
These conditions are stated in terms of the ultragraph
$\G$ and show that, as with graph algebras, the structure of
$C^*(\G)$ is reflected in $\G$.  Using our results from the
previous section we also show that the dichotomy of simple
graph algebras holds for simple ultragraph algebras;
that is, every simple ultragraph algebra is either AF or
purely infinite.

In \S\ref{neither} we use the techniques developed in our
analysis of ideals in \S\ref{simp-sec} to produce an
ultragraph algebra that is neither an Exel-Laca algebra nor
a graph algebra.  This result is important because it shows
that the class of ultragraph algebras is larger than the
graph algebras and the Exel-Laca algebras.  Hence our
results in this paper and the results of \cite{Tom3} are
seen to be more substantial since they hold for
$C^*$-algebras other than just the graph algebras and
Exel-Laca algebras.

We conclude in \S\ref{Cun-Pim} by showing that the
$C^*$-algebras of ultragraphs with no sinks may be
realized as Cuntz-Pimsner algebras.  There is currently much
interest in Cuntz-Pimsner algebras, and since ultragraph
algebras are contained in this class it is possible that
they could serve as interesting examples and perhaps provide
greater insight into the study of general Cuntz-Pimsner
algebras.

%%%%%%%%%%%%%%%%%%%%%%%%%%%%%%%%%%%%%%%%%%%%%%%%%%%%%%%%%%%%%%%
\section{Ultragraph Algebras}
%%%%%%%%%%%%%%%%%%%%%%%%%%%%%%%%%%%%%%%%%%%%%%%%%%%%%%%%%%%%%%%

In this section we review the basic definitions and
properties of ultragraphs and their $C^*$-algebras.  For a
more thorough introduction, we refer the reader to
\cite{Tom3}.

\begin{definition}  An \emph{ultragraph} $\G = (G^0, \G^1,
r, s)$ consists of a countable set of vertices $G^0$, a
countable set of edges $\G^1$, and functions $s : \G^1
\rightarrow G^0$ and $r : \G^1 \rightarrow P(G^0)$, where
$P(G^0)$ denotes the collection of nonempty subsets of
$G^0$. 
\end{definition}

If $\G$ is an ultragraph, then a vertex $v \in
G^0$ is called a \emph{sink} if $| s^{-1}(v) | = 0$ and an
\emph{infinite emitter} if $| s^{-1} (v) | = \infty$.  We
call a vertex a \emph{singular vertex} if it is
either a sink or an infinite emitter.

For an ultragraph $\G = (G^0, \G^1, r, s)$ we let $\G^0$
denote the smallest subcollection of $\Pow(G^0)$ that
contains $\{v \}$ for all $v \in G^0$, contains $r(e)$ for
all $e \in \G^1$, and is closed under finite intersections
and finite unions.  The following lemma gives us another
description of $\G^0$.

\begin{lemma}[\cite{Tom3}, Lemma 2.12 ]  
If $\G := ( G^0,
\G^1,r,s)$ is an ultragraph, then \begin{align*} \G^0 = \{
\bigcap_{e
\in X_1} r(e)
\cup \ldots 
\cup \bigcap_{e \in X_n} r(e) \cup F : & \ \text{$X_1,
\ldots, X_n$ are finite subsets of $\G^1$} \\ & \text{ and $F$
is a finite subset of $G^0$} \}.
\end{align*}  
Furthermore, $F$ may be chosen to
be disjoint from $\bigcap_{e \in X_1} r(e) \cup \ldots 
\cup \bigcap_{e \in X_n} r(e)$.
\label{description}
\end{lemma}

\begin{definition} \label{CK-G-fam} If $\G$ is an ultragraph, a \emph{Cuntz-Krieger $\G$-family} is a collection of
partial isometries $\{ s_e : e \in \G^1 \}$ with mutually
orthogonal ranges and a collection of projections $\{ p_A : A
\in \G^0 \}$ that satisfy
\begin{enumerate}
\item $p_\emptyset = 0$, $p_A p_B = p_{A \cap B}$, and $p_{A
\cup B} = p_A + p_B - p_{A \cap B}$ for all $A,B \in \G^0$
\item $s_e^*s_e = p_{r(e)}$ for all $e \in \G^1$
\item $s_es_e^* \leq p_{s(e)}$ for all $e \in \G^1$
\item $p_v = \sum_{s(e) = v} s_es_e^*$ whenever $0 < |
s^{-1}(v) | < \infty$.
\end{enumerate}

\noindent When $A$ is a singleton set $\{ v \}$, we shall
write $p_v$ in place of $p_{ \{ v \} }$.
\end{definition}

\begin{definition}
If $\G$ is an ultragraph, we let $C^*(\G)$ denote the
$C^*$-algebra generated by a universal Cuntz-Krieger
$\G$-family.  It is proven in \cite[Theorem 2.11]{Tom3} that
$C^*(\G)$ exists.
\end{definition}

\noindent For $n \geq 2$ we define $\G^n := \{ \alpha = \alpha_1 \ldots
\alpha_n : \alpha_i \in \G^1 \text{ and } s(\alpha_{i+1}) \in
r(\alpha_i) \}$ and $\G^* := \bigcup_{n=0}^\infty \G^n$.  The
map $r$ extends naturally to $\G^*$, and we say that $\alpha$
has length $|\alpha| = n$ when $\alpha \in \G^n$.  Note that
the paths of length zero are the elements of $\G^0$, and
when $A \in \G^0$ we define $s(A)=r(A)=A$.

If $\G$ is an ultragraph, then a \emph{loop} is
a path
$\alpha \in \G^*$ with $| \alpha | \geq 1$ and $s(\alpha)
\in r(\alpha)$.  An \emph{exit} for a loop is one of the
following:
\begin{enumerate}
\item an edge $e \in
\G^1$ such that there exists an $i$ for which $s(e) \in
r(\alpha_{i})$ but $e \neq \alpha_{i+1}$
\item a sink $w$ such that $w \in r(\alpha_i)$ for some $i$.
\end{enumerate}

\text{}

\noindent \textbf{Condition~(L):}  Every loop in $\G$ has an
exit; that is, for any loop $\alpha := \alpha_1 \ldots
\alpha_n$ there is either an edge $e \in \G^1$ such that
$s(e) \in r(\alpha_{i})$ and $e \neq
\alpha_{i+1}$ for some $i$, or there is a sink $w$ with $w
\in r(\alpha_i)$ for some $i$.

\text{}

We mention that versions of the Cuntz-Krieger uniqueness
theorem and the gauge-invariant uniqueness theorem have been
proven for ultragraph algebras \cite[Theorem 6.7 and
Theorem 6.8]{Tom3}.

\begin{definition}
If $I$ is a countable set and $A$ is an $I
\times I$ matrix with entries in $\{ 0, 1\}$, then we may
form the ultragraph $\G_A := (G_A^0, \G_A^1, r, s)$
defined by $G_A^0 := \{v_i : i \in I \}$, $\G_A^1 :=  I$,
$s(i) = v_i$ for all $i \in I$, and $r(i)=\{v_j : A_\G (i,j)
= 1 \}$.  
\label{edgeLG}
\end{definition}

Note that the edge matrix of $\G_A$ is $A$.  If $A$ is a
countable $\{0,1\}$-matrix, then it was shown in
\cite[Theorem 4.5]{Tom3} that the Exel-Laca algebra $\OA$
is canonically isomorphic to $C^*(\G)$.

In \cite{EL} Exel and Laca associated a graph
$\Gr(A)$ to $A$ whose vertex matrix is equal to $A$. 
Specifically, one defines the vertices of $\Gr (A)$ to be
$I$, and for each pair of vertices $i,j \in I$ one defines
there to be $A(i,j)$ edges from $i$ to $j$.  We shall see
that the ultragraph $\G_A$ can often tell us more about
the structure of $\OA \cong C^*(\G)$ than the graph $\Gr
(A)$ can.

%%%%%%%%%%%%%%%%%%%%%%%%%%%%%%%%%%%%%%%%%%%%%%%%%%%%%%%%%%%%%%%
\section{Simplicity of Ultragraph Algebras}
\label{simp-sec}
%%%%%%%%%%%%%%%%%%%%%%%%%%%%%%%%%%%%%%%%%%%%%%%%%%%%%%%%%%%%%%%

In \cite[\S4]{BPRS} the ideals of graph algebras were
studied using saturated hereditary subsets of
$G^0$.  Our methods in this section will be similar, except
that we now use saturated hereditary subcollections
of $\G^0$.  Although we could call these subsets of $\G^0$,
we will refer to them as subcollections to emphasize that
their elements are themselves subsets of $G^0$.

\begin{definition}  A subcollection $\Hi \subset \G^0$ is
\emph{hereditary} if 
\begin{enumerate}
\item whenever $e$ is an edge with $\{ s(e) \} \in \Hi$, then
$r(e)
\in \Hi$
\item $A,B \in \Hi$, implies $A \cup B \in \Hi$
\item $A \in \Hi$, $B \in \G^0$, and $B \subseteq A$,
imply that $B \in \Hi$.
\end{enumerate}
\end{definition}

\begin{definition}  A hereditary subcollection $\Hi \subset
\G^0$ is \emph{saturated} if for any $v \in G^0$ with $0
< |s^{-1}(v)| < \infty$ we have that $$\{ r(e) : e \in \G^1
\text{ and } s(e)=v \} \subseteq \Hi \text{
\quad implies \quad } \{v \} \in \Hi.$$
\end{definition}

\noindent The \emph{saturation} of a hereditary collection
$\Hi$ is the smallest saturated subcollection
$\overline{\Hi}$ of
$\G^0$ containing $\Hi$; the saturation $\overline{\Hi}$ is
itself hereditary.

\begin{remark}  Note that if $\Hi \subseteq \G^0$ is a
hereditary subcollection with $\{ v \} \in \Hi$ for all $v
\in G^0$, then $\Hi = \G^0$.  This is because having
$\Hi$ hereditary implies that $\Hi$ contains $r(e)$ for
all $e \in \G^1$, and since $\Hi$ is closed under
finite unions and intersections Lemma~\ref{description} then
implies $\Hi = \G^0$.
\end{remark}

\begin{lemma}  Let $\G$ be an ultragraph and let $I$ be an
ideal in $C^*(\G)$.  Then $\Hi := \{ A \in \G^0 : p_A \in I
\}$ is a saturated hereditary subcollection of $\G^0$.
\label{sat-her}
\end{lemma}

\begin{proof} Suppose $\{ s(e) \} \in \Hi$.  Then
$$p_{s(e)} \in I \Longrightarrow s_e = p_{s(e)}s_e \in I
\Longrightarrow p_{r(e)} = s_e^*s_e \in I \Longrightarrow
r(e) \in \Hi.$$ Also, if $A,B \in \Hi$,
then $$p_A, p_B \in I \Longrightarrow p_{A \cup B} = p_A +
p_B - p_Ap_B \in I \Longrightarrow A \cup B \in \Hi.$$ 
Finally, if $A \in H$, $B \in \G^0$, and
$B \subseteq A$, then $$p_A \in I \Longrightarrow p_B =
p_Bp_A \in I \Longrightarrow B \in H$$ so $H$ is hereditary. 

Furthermore, if $0 < |s^{-1}(v)| < \infty$
and $\{r(e) : \text{$e \in \G^1$ and $s(e)=v$} \} \subseteq
H$, then $\{ s_e : \text{$e \in \G^1$ and $s(e)=v$} \}
\subseteq I$ and $p_v = \sum_{s(e)=v} s_es_e^* \in I$ which
implies that $\{ v \} \in \Hi$.  Thus
$H$ is saturated.
\end{proof}

\noindent  For a hereditary subcollection $\Hi
\subseteq \G^0$ let $I_{\Hi}$ denote the ideal in $C^*(\G)$
generated by $\{ p_A : A \in \Hi \}$.

\begin{lemma}  Let $\G$ be an ultragraph and let $\Hi$ be a
hereditary subcollection of $\G^0$.  Then $$I_{\Hi} =
\overline{\mathrm{span}} \{ s_\alpha p_A s_\beta^* :
\alpha,\beta \in \G^* \text{ and } A \in
\overline{\Hi}  \}.$$  In particular, $I_{\Hi} =
I_{\overline{\Hi}}$ and $I_{\Hi}$ is gauge invariant.
\label{IH}
\end{lemma}

\begin{proof} Note that $\{ A \in \G^0 : p_A \in I_H \}$ is
a saturated set containing $\Hi$ and therefore contains
$\overline{\Hi}$.  Thus $J := \overline{\mathrm{span}} \{ s_\alpha p_A s_\beta^* :
\alpha,\beta \in \G^* \text{ and } A \in
\overline{\Hi}  \}$ is contained in $I_\Hi$.  For inclusion
in the other direction, notice that any nonzero product of
the form $s_\alpha p_A s_\beta^* s_\gamma p_B s_\delta^*$
collapses to another of the form $s_\mu p_C s_\nu^*$ and
from an examination of the various possibilities and the
hereditary property of $\overline{\Hi}$ we deduce that $J$
is an ideal.  Since $J$ contains the generators of $I_\Hi$,
it follows that $J := I_{\Hi}$.  The last two remarks follow
easily.
\end{proof}

\begin{lemma}  Let $\G$ be an ultragraph for which
$C^*(\G)$ is simple.  If $\Hi$ is a saturated hereditary
subcollection of $\G^0$ and $K := \{ v \in G^0 : \{v\} \in
\Hi \}$, then for any $e \in \G^1$ we have that $r(e)
\subseteq K$ implies that $r(e) \in \Hi$.
\label{simplerange}
\end{lemma}

\begin{proof}  If $\Hi$ is empty the claim holds
vacuously.  If $\Hi \neq \emptyset$, then since
$C^*(\G)$ is simple we know that $I_\Hi = C^*(\G)$ and thus
$p_{r(e)} \in I_\Hi$.  By Lemma~\ref{IH} there exist
$\lambda_k \in \C$, $\alpha_k,\beta_k \in \G^*$, and
$A_k,B_k \in \Hi$ for $1 \leq k \leq n$ such that
$$\| p_{r(e)} - \sum_{k=1}^n \lambda_k s_{\alpha_k}
p_{A_k} s_{\beta_k}^*  \| < 1.$$
Furthermore, since 
$$\| p_{r(e)} \Big( p_{r(e)} - \sum_{k=1}^n \lambda_k
s_{\alpha_k} p_{A_k} s_{\beta_k}^* \Big) \| \leq \| p_{r(e)}
- \sum_{k=1}^n \lambda_k s_{\alpha_k} p_{A_k} s_{\beta_k}^* 
\|$$ we may assume that $s(\alpha_k) \in r(e)$ when
$| \alpha_k | \geq 1$ and $s(\alpha_k) \subseteq r(e)$ when
$ | \alpha_k | = 0$.  (We remind the reader that if $|\alpha|
= 0$, then $\alpha = A$ for some $A \in \G^0$ and
$s(\alpha):= A$.)

Now define $B := \bigcup_{k=1}^n s(\alpha_k)$.  Since $B
\subseteq r(e)$ we see that $q:= p_{r(e)} - p_B$ is a
projection.  Furthermore, 
$$ \| q \| = \| q \Big( p_{r(e)} - \sum_{k=1}^n
\lambda_k s_{\alpha_k} p_{A_k} s_{\beta_k}^* \Big) \| 
\leq \|  p_{r(e)} - \sum_{k=1}^n
\lambda_k s_{\alpha_k} p_{A_k} s_{\beta_k}^*  \| < 1. $$
and since $q$ is a projection this implies that $q=0$. 
Therefore $p_{r(e)} = p_B$ and $r(e) = B =
\bigcup_{k=1}^n s(\alpha_k) \in \Hi$. 
\end{proof}

\begin{lemma}  Let $\G$ be an ultragraph for which
$C^*(\G)$ is simple.  If $\Hi$ is a saturated
hereditary subcollection of $\G^0$, then either $\Hi = \G^0$
or $\Hi = \emptyset$.
\label{her-simp}
\end{lemma}

\begin{proof}  Let $\G = (G^0,\G^1,r,s)$.  Set $K := \{ w
\in G^0 : \{w\} \in \Hi \}$ and $S := G^0 \backslash K$. 
We define an ultragraph $\F = (F^0, \F^1, r_{\F},
s_{\F})$ as follows:
\begin{align*}
F^0 & := S & & s_{\F}(e)  := s(e) \\
\F^1 & := \{ e \in \G^1 : r(e) \cap S \neq \emptyset \} & &
r_{\F}(e)  := r(e) \cap S 
\end{align*} Note that if $e \in \F^1$, then $r(e) \cap S
\neq \emptyset$ so $r(e) \notin \Hi$ and since $\Hi$ is
hereditary it follows that $\{ s(e)\}  \notin \Hi$ and $s(e)
\in S$.  Thus $s_{\F}$ is well-defined. 

Let $\{ s_e,p_A \}$ be the canonical Cuntz-Krieger
$\F$-family in $C^*(\F)$.  For each $e \in \G^1$ and $A
\in \G^0$ define
$$t_e := \begin{cases} s_e & \text{if $e \in
\F^1$} \\ 0 & \text{otherwise} \end{cases} \quad \text{ and
} \quad q_A := p_{A \cap S}.$$ Note that if $A \in \G^0$,
then by Lemma~\ref{description} $$A = \bigcap_{e \in X_1}
r(e) \cup \ldots \cup \bigcap_{e \in X_n} r(e) \cup F$$ for
some finite subsets $X_1, \ldots, X_n \subseteq \G^1$ and
some finite subset $F \subseteq G^0$.  Thus if $$Y_i
:= \begin{cases} X_i & \text{if $X_i \subseteq
\F^1$} \\ \emptyset & \text{otherwise} \end{cases}$$
then we see that 
\begin{align*} A \cap S & = \bigcap_{e \in X_1}
(r(e) \cap S) \cup \ldots \cup \bigcap_{e \in X_n} (r(e)
\cap S) \cup (F \cap S) \\
& = \bigcap_{e \in Y_1}
r_{\F}(e) \cup \ldots \cup \bigcap_{e \in Y_n} r_{\F}(e) \cup
(F \cap S)
\end{align*}
which is in $\F^0$.  Hence $q_A$ is well-defined.

We shall now show that $\{t_e,q_A\}$ is a Cuntz-Krieger
$\G$-family.  Clearly, the $t_e$'s have mutually orthogonal
ranges since the $s_e$'s do.  Thus we simply need to verify
the four properties of Definition~\ref{CK-G-fam}.
\begin{enumerate}
\item We have that $q_\emptyset = p_\emptyset =
0$, $q_Aq_B = p_{A \cap S} p_{B \cap S} = p_{(A \cap B) \cap
S} = q_{A \cap B}$, and $q_{A \cup B}  = p_{(A \cup B) \cap
S} = p_{(A \cap S) \cup (B \cap S)} = p_{(A \cap S)} + p_{(B
\cap S)} - p_{(A \cap S) \cap (B \cap S)} = q_A + q_B
-p_{(A\cap B) \cap S} = q_A + q_B -q_{A \cap B}$.
\item If $e \in \F^1$, then $t_e^*t_e = s_e^*s_e =
p_{r_{\F}(e)} = p_{r(e) \cap S} = q_{r(e)}$.  On the other
hand, if $e \notin \F^1$, then $r(e) \cap S = \emptyset$ so
$q_{r(e)} = 0 = t_e^* t_e$.
\item  If $e \in \F^1$, then $s(e) \in S$ so $t_et_e^* =
s_es_e^* \leq p_{s_{\F}(e)} = q_{s(e)}$.  On the other hand,
if $e \notin \F^1$, then $t_et_e^* = 0 \leq q_{s(e)}$.
\item  Let $v \in G^0$ and $0 < | s^{-1}(v) | <
\infty$.  If $v \notin S$, then $\{ v \} \in \Hi$ and $r(e)
\in \Hi$ so $r(e) \cap S = \emptyset$ and 
$$\sum_{\{ e \in \G^1 : s(e) = v \} } t_et_e^* = 0 =
q_{r(e)}.$$  If $v \in S$, then since $s_{\F}^{-1}(v)
\subseteq s^{-1}(v)$ we have that $| s_{\F}^{-1}(v) | <
\infty$.  Also note that $r(e) \cap S = \emptyset$ implies
$r(e) \in \Hi$ by Lemma~\ref{simplerange}.  Thus the fact
that $\{v\} \notin \Hi$ and the fact that $\Hi$ is saturated
imply that there is at least one edge $e$ with
$r(e) \cap S \neq \emptyset$.  Hence $0 < | s_{\F}^{-1}(v)
|$.  Thus
\begin{align*} \sum_{\{ e \in \G^1 : s(e) = v \} } t_et_e^*
& = \sum_{\{ e \in \F^1 : s(e) = v \} } t_et_e^* + \sum_{\{ e
\in (\G^1 \backslash \F^1) : s(e) = v \} } t_et_e^* \\
& = \sum_{\{ e \in \F^1 : s_{\F}(e) = v \} } s_es_e^* +
0 = p_{s_{\F}(e)} = q_{s(e)}.
\end{align*}
\end{enumerate}
Now since $\{t_e, q_v\}$ is a Cuntz-Krieger $\G$-family with
$q_v = 0$ if and only if $v \notin S$, the universal property
gives a homomorphism $\phi : C^*(\G) \rightarrow C^*(\F)$
whose kernel contains only those projections corresponding
to vertices that are not in $S$.  Since $C^*(\G)$ is simple,
the kernel of $\phi$ is either $C^*(\G)$ or $\{ 0 \}$.  Thus
$S$ is either $\emptyset$ or $G^0$, and $K$ is either $G^0$
or $\emptyset$.  Since $\Hi$ is a saturated hereditary
subset, this implies that either $\Hi = \emptyset$ or $\Hi =
\G^0$.
\end{proof}

The following proof is modeled after that of \cite[Theorem
4.1(c)]{BPRS}.

\begin{lemma} Let $\G$ be an ultragraph.  If $\G$ has a
loop with no exits, then $C^*(\G)$ contains an ideal Morita
equivalent to $C( \T )$.
\label{simple-exits-nosinks}
\end{lemma}

\begin{proof}  Let $C^*(\G) = C^*(\{s_e,p_A \})$ and
$\alpha = \alpha_1 \ldots \alpha_n$ be a loop in $\G$ with no
exits.  This implies that
$r(\alpha_i) = \{ s(\alpha_{i+1}) \}$ for $1
\leq i < n$, and $r(\alpha_n) = \{ s(\alpha_1) \}$.  In
particular, the $\alpha_i$'s have ranges that are singleton
sets.  Define $X := \{ s(\alpha_i) \}_{i=1}^n$ and $q_X :=
\sum_{v \in X} p_v$.  If $\Hi$ equals the (finite) collection
of all subsets of $X$, then $\Hi$ is a hereditary subset of
$\G^0$.  We shall show that $I_\Hi = I_{\overline{\Hi}}$ is
Morita equivalent to $C(\T)$.

Define $G^1 := \{ \alpha_i \}_{i=1}^n$ and let $G$ be the
graph $G := (X, G^1, r, s )$.  We claim that
$q_X I_{\overline{\Hi}} q_X$ is generated by the
Cuntz-Krieger $G$-family $\{ s_e, p_v : e \in G^1, v \in X
\}$.  Certainly this family lies in the corner.  On the
other hand, if $\alpha, \beta \in \G^*$ and $A \in \G^0$,
then $q_X s_\alpha p_A s_\beta^* q_X = 0$ unless both
$\alpha$ and $\beta$ have sources in $X$.  Thus the claim is
verified and the gauge-invariant uniqueness theorem for
$C^*$-algebras of graphs \cite[Theorem 2.1]{BPRS} implies
that $q_X I_{\overline{\Hi}} q_X \cong C(G)$.  To see that
this is a full corner of $I_\Hi$, suppose that $J$ is an
ideal in $I_\Hi$ containing $q_X I_{\overline{\Hi}} q_X$. 
Then $J$ is an ideal of $C^*(\G)$ and Lemma~\ref{sat-her}
implies that $\{ A \in \G^0 : p_A
\in J \}$ is a saturated hereditary subcollection containing
$\{ \{v \} : v \in X \}$ and hence containing $\Hi$.  But
this implies that $J$ contains the generators of
$I_{\overline{\Hi}}$ and hence is all of
$I_{\overline{\Hi}}$.  

Therefore, $I_{\overline{\Hi}}$ is Morita equivalent to
$C^*(G)$.  Since $G$ is a loop of length $n$, we see
from \cite[Theorem 2.4]{KPR} that $C^*(G) \cong C( \T)
\otimes M_n(\C)$ which is Morita equivalent to $C(\T)$.
\end{proof}

\begin{lemma} Let $\G$ be an ultragraph such that
$C^*(\G)$ is simple.  Then every loop in $\G$ has an
exit.
\label{simple-exits}
\end{lemma}

\begin{proof}  If $\G$ contained a loop with no exits, then
Lemma~\ref{simple-exits-nosinks} would imply that $C^*(\G)$
contains an ideal Morita equivalent to $C(\T)$.  Hence
$C^*(\G)$ could not be simple. 
\end{proof}

The following is a generalization of \cite[Theorem 8]{Szy}
\begin{theorem}  If $\G$ is an ultragraph,
then $C^*(\G)$ is simple if and only if $\G$ satisfies:
\begin{enumerate}
\item every loop in $\G$ has an exit
\item the only saturated hereditary subcollections of $\G^0$
are $\G^0$ and $\emptyset$.
\end{enumerate}
\label{Woj-simplicity}
\end{theorem}

\begin{proof}  Suppose that $C^*(\G)$ is simple.  Then
Lemma~\ref{simple-exits} implies that every loop in $\G$
must have an exit.  Furthermore, if $\Hi$ is a saturated
hereditary subcollection of $\G^0$, then it follows from
Lemma~\ref{her-simp} that either $\Hi = \G^0$ or $\Hi =
\emptyset$.

Conversely, suppose that $\G$ satisfies the two properties
above.  If $I$ is an ideal in $C^*(\G)$, then
Lemma~\ref{sat-her} tells us that $\Hi := \{ A \in \G^0 :
p_A \in I \}$ is a saturated hereditary subcollection of
$\G^0$.  Hence $\Hi$ equals either $\G^0$ or $\emptyset$. 
If $\Hi = \G^0$, then clearly $I = C^*(\G)$.  On the other
hand, if $\Hi = \emptyset$, then since every loop in
$\G$ has an exit we may use the Cuntz-Krieger Uniqueness
Theorem to conclude that the projection $\pi : C^*(\G)
\rightarrow C^*(\G) / I$ is injective, and thus $I = \{ 0
\}$.
\end{proof}

Recall from \cite[Corollary~2.15]{DT} and
\cite[Theorem~4]{Pat} that if
$G$ is a graph, then $C^*(G)$ is simple if and only if
every loop in $G$ has an exit, $G$ is cofinal, and $G^0 \geq
\{v \}$ for every singular vertex $v \in G^0$.  We shall use
the previous theorem to obtain an analogous characterization
for ultragraph algebras.  Recall that infinite emitters
in a graph correspond to infinite sets of the form $r(e) \in
\G^0$.  In fact, if $G$ is a graph with vertex matrix $A$,
then in the ultragraph $\G_A$ the set $r(e)$ is finite
for all $e \in \G_A^1$ if and only if $G$ has no infinite
emitters.

We first extend the notions of $\geq$ and cofinality to
ultragraphs.  If $\G$ is an ultragraph and
$v,w \in G^0$, we write $w \geq v$ to mean that there exists
a path $\alpha \in \G^*$ with $s(\alpha) = w$ and $v \in
r(\alpha)$.  Also, we write $G^0 \geq \{v \}$ to mean that
$w \geq v$ for all $w \in G^0$.  We say that $\G$ is
\emph{cofinal} if for every infinite path $\alpha := e_1 e_2
\ldots$ and every vertex $v
\in G^0$ there exists an $i \in \N$ such that $v
\geq s(e_i)$.

In addition, we need a new notion of reachability.  If $v
\in G^0$ and $A \subseteq G^0$, then we write $v \rightarrow
A$ to mean that there exist a finite number of
paths $\alpha_1, \ldots \alpha_n \in
\G^*$ such that $s(\alpha_i) = v$ for all $1 \leq i \leq n$
and $A \subseteq \bigcup_{i=1}^n r(\alpha_i)$.  Note that if
$A = \{ w \}$, then $v \rightarrow \{w \}$ if and only if $v
\geq w$.

\begin{theorem}
If $\G$ is an ultragraph,
then $C^*(\G)$ is simple if and only if $\G$ satisfies:
\begin{enumerate}
\item every loop in $\G$ has an exit
\item $\G$ is cofinal
\item $G^0 \geq \{ v\}$ for every singular vertex $v \in G^0$
\item If $e \in \G^1$ is an edge for which the set
$r(e)$ is infinite, then for every $w \in G^0$ there exists a
set $A_w \subseteq r(e)$ for which $r(e) \backslash A_w$ is
finite and $v \rightarrow A_w$.
\end{enumerate}
\label{DT-simplicity}
\end{theorem}

\noindent In order to prove this result we need a lemma.

\begin{lemma}  Let $\G$ be an ultragraph and let $\Hi
\subseteq \G^0$ be a hereditary subset.  Set $\Hi_0 := \Hi$
and for $n \in \N$ define $$\Hi_{n+1} := \{ A \cup F : A \in
\Hi_{n} \text{ and $F$ is a finite subset of $S_n$} \}$$
where $S_n := \{ w \in G^0 : 0 < | s^{-1}(w) | < \infty
\text{ and } \{ r(e) : s(e) = w \} \subseteq \Hi_n \}$.
Then $$\overline{\Hi} = \bigcup_{i=0}^\infty
\Hi_i$$ and every $X \in \overline{\Hi}$
has the form $X = A \cup F$ for some $A \in \Hi$ and some
finite set $F \subseteq \bigcup_{i=1}^\infty S_i$.
\label{induct-sat}
\end{lemma}

\begin{proof}  To see that $\bigcup_{i=0}^\infty \Hi_i
\subseteq \overline{\Hi}$, first note that $\Hi_0 \subseteq
\overline{\Hi}$.  Also whenever $\Hi_n \subseteq
\overline{\Hi}$, then because $\overline{\Hi}$ is saturated
we have that $F \in \overline{\Hi}$ for any finite subset $F
\subseteq S_n$, and hence $\Hi_{n+1} \subset \overline{\Hi}$.
Thus by induction we have $\Hi_n \subseteq \overline{\Hi}$
for all $n$.

To see that $\overline{\Hi} \subseteq \bigcup_{i=0}^\infty
\Hi_i$ we shall show that $\bigcup_{i=0}^\infty \Hi_i$ is a
saturated hereditary subcollection.  We shall begin by
proving inductively that each $\Hi_i$ is hereditary.  For
the base case, we have by hypothesis that $\Hi_0 := \Hi$ is
hereditary.  Now assume that $\Hi_{n}$ is hereditary and
consider $\Hi_{n+1}$.  If $\{s(e) \} \in
\Hi_{n+1}$, then by the definition of $\Hi_{n+1}$ either $\{
s(e) \} \in \Hi_{n}$ or $ s(e) \in S_{n}$.  In either
case, $r(e) \in \Hi_{n} \subseteq \Hi_{n+1}$. To see that
$\Hi_{n+1}$ is closed under subsets, Let $A_1 \cup F_1$ and
$A_2 \cup F_2$ be typical elements of $\Hi_{n+1}$.  Then
$(A_1 \cup F_1) \cup (A_2 \cup F_2) = (A_1 \cup A_2) \cup
(F_1 \cup F_2)$ which is in $\Hi_{n+1}$ because $\Hi_{n}$
is closed under unions.  Finally, suppose that $A \cup F$ is
a typical element of $\Hi_{n+1}$ and that $B \in \G^0$ with
$B \subseteq A \cup F$.  Then $A \cap B \subseteq A$ and
since $\Hi_{n}$ is hereditary $A \cap B \in \Hi_{n}
\subseteq \Hi_{n+1}$.  Also,
$B \cap F \subseteq F$ so $B \cap F$ is a finite subset of
$S_n$.  Thus $B = (B \cap A) \cup (B \cap F) \in \Hi_{n+1}$
and $\Hi_{n+1}$ is hereditary.

Since $\bigcup_{i=0}^\infty \Hi_i$ is the union of
hereditary sets, it follows that $\bigcup_{i=0}^\infty
\Hi_i$ itself is hereditary.  To see that
$\bigcup_{i=0}^\infty \Hi_i$ is also saturated, let $v \in
G^0$ be a vertex with $0 < | s^{-1}(v) | < \infty$ and $\{
r(e) : s(e) = v \} \subseteq \bigcup_{i=0}^\infty \Hi_i$. 
Since $\Hi_i \subseteq \Hi_{i+1}$ and since there are only
finitely many edges with source $v$, we see that there
exists $n \in \N$ such that $\{ r(e) : s(e) = v \} \subseteq
\Hi_n$.  Thus $v \in S_n$ and $\{v\} \in \Hi_{n+1}
\subseteq \bigcup_{i=0}^\infty \Hi_i$.  Hence
$\bigcup_{i=0}^\infty \Hi_i$ is saturated.  

Therefore $\overline{\Hi} = \bigcup_{i=0}^\infty \Hi_i$ and
to prove the claim, we let $X \in \bigcup_{i=0}^\infty
\Hi_i$.  Then $X \in \Hi_n$ for some $n \in \N$ and $X =
A_{n-1} \cup F_{n-1}$ for some $A_{n-1} \in \Hi_{n-1}$ and
some finite subset $F_{n-1} \subseteq S_{n-1}$.  Similarly,
$A_{n-1} = A_{n-2} \cup F_{n-2}$ for some $A_{n-2} \in
\Hi_{n-2}$ and some finite subset $F_{n-2}
\subseteq S_{n-2}$.  Continuing inductively we see that $A =
A_0 \cup (F_{n-1} \cup \ldots \cup F_1)$ where the $F_i$'s
are all finite sets.
\end{proof}

\noindent \emph{Proof of Necessity in Theorem
\ref{DT-simplicity}.}  Suppose that $C^*(\G)$ is simple.  By
Theorem~\ref{Woj-simplicity} we see that every loop in
$\G$ has an exit and the only saturated hereditary subsets of
$\G^0$ are $\G^0$ and $\emptyset$.  

Let $\alpha = e_1 e_2 \ldots $ be an infinite path
and set $K := \{ w \in G^0 : w \ngeq s(e_i) \text{
for all $i$} \}$.  Also define $\Hi := \{ A \in \G^0 : A
\subseteq K \}$.  Then one can verify that $\Hi$ is a
saturated hereditary subcollection of $\G^0$.  Since $\{
s(e_1) \} \notin \Hi$ we see that $\Hi$ is not all of
$\G^0$.  Thus $\Hi = \emptyset$ and $\G$ is cofinal.

Let $v \in G^0$ be a singular vertex.  Fix any vertex $w
\in G^0$ and define $K := \{ x \in G^0 : w \geq x \}$.  Also
let $\Hi : = \{ A \in \G^0 : A \subseteq K \}$.  Then
$\Hi$ is a hereditary subcollection of $\G^0$.  If
$\overline{\Hi}$ is the saturation of $\Hi$, then
$\overline{\Hi}$ is nonempty because $\{w \} \in \Hi$. 
Hence $\overline{\Hi} = \G^0$.  Now, using the notation of
Lemma~\ref{induct-sat}, we see that $v \notin S_i$ for all
$i$ because $v$ is a singular vertex.  Therefore it follows
from Lemma~\ref{induct-sat} that $\{ v \} \in \overline{\Hi}$
implies that $\{ v \} \in \Hi$.  Thus $v \in K$ and $w \geq
v$.  Hence $G^0 \geq \{ v \}$.

Let $e \in \G^1$ be an edge such that $r(e)$ is an infinite
set.  Fix $w \in G^0$ and set $\Hi := \{ A \in \G^0 : w
\rightarrow A \}$.  To see that $\Hi$ is hereditary suppose
that $\{ s(f) \} \in \Hi$.  Then $w \rightarrow \{s(f) \}$
and hence $v \geq s(f)$.  Thus there exists a path $\beta$
with $s(\beta) = w$ and $s(f) \in r(\beta)$.  But then we see
that $w \rightarrow r(f)$ via the path $\beta f$. 
Additionally, it is easy to see that $\Hi$ is closed under
unions and subsets.  Since $\{w\} \in \Hi$, it follows that
$\Hi$ is nonempty, and hence $\overline{\Hi} = \G^0$.  Thus
$r(e) \in \overline{\Hi}$.  By Lemma~\ref{induct-sat} it
follows that $r(e) = A_w \cup F$ for some $A_w \in \Hi$ and
some finite set $F$.  But then $w \rightarrow A_w$ and
$r(e) \backslash A_w$ is finite. 
\hfill $\square$

\text{ }

\noindent \emph{Proof of Sufficiency in
Theorem~\ref{DT-simplicity}.}  Suppose that $\G$ satisfies
the four conditions stated in Theorem~\ref{DT-simplicity}. 
In light of Theorem~\ref{Woj-simplicity} it suffices to show
that the only saturated hereditary subcollections of $\G^0$
are $\emptyset$ and $\G^0$.

Let $\Hi$ be a nonempty saturated hereditary subcollection
of $\G^0$.  We shall show
that for every $w \in G^0$ with $\{ w \} \notin \Hi$ there
exists an edge $e \in \G^1$ such that
$s(e) = w$ and $r(e)$ contains a vertex $w'$ for which $\{
w' \} \notin \Hi$.

If $\{ w \} \notin \Hi$, then since $G^0 \geq \{ v \}$ for
every singular vertex $v$ it follows that $w$ is not a
singular vertex.   Therefore, since $\Hi$ is saturated, there
exists an edge $e$ such that $s(e) = w$ and $r(e)
\notin \Hi$.  If $r(e)$ is finite, then because $\Hi$ is
closed under unions, there must exist a vertex $w' \in r(e)$
such that $\{ w' \} \notin \Hi$.  If $r(e)$ is infinite, then
choose some $x \in G^0$ for which $\{x \} \in \Hi$.  Then
there exists $A_x \subseteq r(e)$ such that $w \rightarrow
A_x$ and $r(e) \backslash A_x$ is a finite set.  Let
$\alpha_1, \ldots, \alpha_n$ be paths with
$s(\alpha_i) = x$ and $A_x \subseteq \bigcup_{i=1}^n
r(\alpha_i)$.  Since $\Hi$ is hereditary, it follows that
$\bigcup_{i=1}^n r(\alpha_i) \in \Hi$.  Now we must have that
one of the vertices in $r(e) \backslash A_x$ is not in
$\Hi$.  For otherwise, $r(e) \backslash A_x \in \Hi$ and
$\bigcup_{i=1}^n r(\alpha_i) \cup (r(e) \backslash A_x)$ is
an element in $\Hi$ containing $r(e)$ which contradicts the
fact that $r(e) \notin \Hi$.

Now suppose that there exists $w_1 \in G^0$ such that $\{ w_1
\} \notin \Hi$.  From the argument in the preceding
paragraph there exists an edge $e_1$ and a vertex $w_2$
such that $s(e_1) = w_1$, $w_2 \in r(e_1)$, and $\{ w_2 \}
\notin \Hi$.  Continuing inductively, we create an
infinite path $e_1e_2e_3 \ldots$ with $\{ s(e_i) \} \notin
\Hi$ for all $i$.  But this contradicts the cofinality of
$\G$.  Hence $\Hi$ must be all of $\G^0$.
\hfill $\square$

\begin{remark}  Note that if $\G$ is an ultragraph with
two sinks $v_1$ and $v_2$, then $v_1 \ngeq v_2$ and hence
$G^0 \ngeq v_2$.  Therefore if $\G$ is an ultragraph with
two or more sinks, then $C^*(\G)$ is not simple.  In
addition, if $\G$ has exactly one sink and $C^*(\G)$ is
simple, then $\G$ contains no infinite paths because the
sink is unable to reach any infinite path.
\end{remark}

\begin{example} 
\label{simp-ex} Let $\G$ be the ultragraph
$$
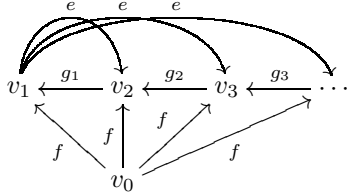
\xymatrix{
v_1 \ar@(u,u)[r]^{e} \ar@(u,u)[rr]^{e} 
\ar@(u,u)[rrr]^{e} & v_2
\ar[l]_{g_1} & v_3 \ar[l]_{g_2} & \cdots \ar[l]_{g_3} \\
 & v_0 \ar[ul]^f \ar[u]^{f} \ar[ur]^{f}
\ar[urr]_{f}  & & \\ }$$
Since any loop in $\G$ must contain the edge $e$, we see
that every loop has an exit.  In addition, let
$\Hi$ be a saturated hereditary subcollection of $\G^0$.  If
$\Hi$ is nonempty, then there is some singleton set $\{ v \}
\in \Hi$ for $v \in G^0$.  Because $\Hi$ is hereditary and
because each $v_i$ can be reached from any other vertex, we
see that $\Hi$ must contain $\{ v_i \}$ for $1 \leq i <
\infty$.  Since $s(e) = v_1$ and $\{v_1 \} \in \Hi$ it
follows that $r(e) = \{ v_2, v_3, v_4, \ldots \} \in \Hi$. 
Hence $r(f) = \{v_1, v_2, v_3, \ldots \} = \{v_1 \} \cup
r(e) \in \Hi$.  Since $\Hi$ is saturated we also have that
$\{ v_0 \} \in \Hi$.  Thus $\Hi$ contains $r(e)$, $r(f)$, and
$\{v\}$ for all $v \in G^0$.  Consequently $\Hi = \G^0$.  It
follows from Theorem~\ref{Woj-simplicity} that $C^*(\G)$ is
simple.
\end{example}

\begin{example}
Consider the infinite matrix
$$
A = \left(\begin{smallmatrix}
1 & 1 & 1 & 1 & 1 & 1 & \\
0 & 0 & 1 & 1 & 1 & 1 & \\
0 & 1 & 0 & 0 & 0 & 0 & \cdots \\ 
0 & 0 & 1 & 0 & 0 & 0 & \\ 
0 & 0 & 0 & 1 & 0 & 0 & \\ 
0 & 0 & 0 & 0 & 1 & 0 & \\
  &   & \vdots  &  & & & \ddots \\ 
\end{smallmatrix}\right)
$$

\noindent The graph $\Gr (A)$ associated to $A$ is

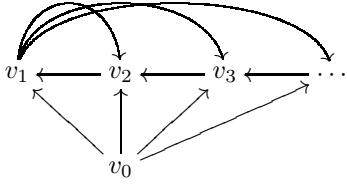
\vspace{.3in}

$$\xymatrix{
v_1 \ar@(u,u)[r] \ar@(u,u)[rr] 
\ar@(u,u)[rrr] & v_2
\ar[l] & v_3 \ar[l] & \cdots \ar[l] \\
 & v_0 \ar[ul] \ar[u] \ar[ur]
\ar[urr]  & & \\ }$$
Since $v_0$ is an infinite emitter, we see that
$H := \{v_1, v_2, v_3, \ldots \}$ is a saturated
hereditary subset of $G^0$ that gives rise to a nontrivial
ideal.  Therefore $C^*(\Gr(A))$ is not simple.

However, the ultragraph $\G_A$ associated to $A$ is the
ultragraph shown in Example~\ref{simp-ex}, and as we saw
there $C^*(\G_A)$ is simple.  Note that $\G_A$ has no
infinite emitters or sinks, and in fact, $| s^{-1} (v) | =1$
for all $v \in G^0$.  It then follows from
\cite[Theorem~4.5]{Tom3} that $C^*(\G_A) \cong \OA$.  Thus
$\OA$ is simple.

This example shows that the ultragraph $\G_A$ is a better
tool for studying $\OA$ than the graph $\Gr(A)$.  As we
saw, the fact that $\OA$ is simple is reflected in the
ultragraph $\G_A$, but it is difficult to see from the
graph $\Gr (A)$. In addition,
$C^*(\Gr(A))$ is very different from $\OA$ whereas
$C^*(\G_A) \cong \OA$. 
\end{example}

%%%%%%%%%%%%%%%%%%%%%%%%%%%%%%%%%%%%%%%%%%%%%%%%%%%%%%%%%%%%%%%
\section{AF and purely infinite ultragraph algebras}
\label{AF-pi}
%%%%%%%%%%%%%%%%%%%%%%%%%%%%%%%%%%%%%%%%%%%%%%%%%%%%%%%%%%%%%%%

\begin{theorem}  Let $\G$ be an ultragraph.  Then
$C^*(\G)$ is an AF-algebra if and only if $\G$ has no
loops.  
\label{AF}
\end{theorem}

\begin{proof}
Let $\F$ be a desingularization of $\G$ \cite[Definition
6.3]{Tom3}.  Then since the class of AF-algebras is closed
under stable isomorphism
\cite[Theorem 9.4]{Eff}, and since $\F$ has loops if and only
if $\G$ has loops, we see that it suffices to prove the
claim for ultragraphs with no singular vertices.

Suppose $\G$ has no singular vertices.  If $\G$ has no
loops, then write $\G^1 := \bigcup_{n=1}^\infty F_n$ as the
increasing union of finite subsets $F_n$, and let $B_n$ be
the $C^*$-subalgebra of $C^*(\G)$ generated by $\{s_e : e
\in F_n \}$.  By \cite[Corollary~5.4]{Tom3} there are
isomorphisms $\phi_n : C^*(G_{F_n}) \rightarrow B_n$.  Since
$\G$ has no loops, it follows from \cite[Lemma~5.6]{Tom3}
that each $G_{F_n}$ has no loops.  Since $G_{F_n}$ is a
finite graph with no loops, $C^*(G_{F_n}) \cong B_n$ is a
finite-dimensional $C^*$-algebra \cite[Corollary
2.3]{KPR}.  Because $\G$ has no singular vertices, the
$s_e$'s are dense in $C^*(\G)$ and $C^*(\G) =
\overline{\bigcup_{n=1}^\infty B_n}$.  Thus $C^*(\G)$ is the
direct limit of finite-dimensional $C^*$-algebras, and
consequently $C^*(\G)$ is an AF-algebra.

Conversely, suppose that $\G$ has a loop $\alpha := \alpha_1
\ldots \alpha_n$.
  
\noindent \textsc{Case I:} $\alpha$ has an exit. 

Because $\G$ has no sinks, we may assume without loss of
generality that there exists an edge $f \in \G^1$ with $f
\neq \alpha_1$ and $s(f) \in r(\alpha_n)$.  Now $$p_{r(e)} =
s_\alpha^*s_\alpha \sim s_\alpha s_\alpha^* \leq
s_{\alpha_1} s_{\alpha_1}^* < s_{\alpha_1} s_{\alpha_1}^* +
s_fs_f^* \leq p_{s(\alpha)} \leq p_{r(\alpha)}$$ and so
$p_{r(\alpha)}$ is an infinite projection.  Since a
projection in an AF-algebra is equivalent to one in a
finite-dimensional subalgebra it cannot be infinite.  Hence
$C^*(\G)$ is not AF.

\noindent \textsc{Case II:} $\alpha$ has no exits. 

Since $\G$ has no sinks, it follows from
Lemma~\ref{simple-exits-nosinks} that $C^*(\G)$ contains an
ideal Morita equivalent to $C(\T)$ which is not AF.  Hence
$C^*(\G)$ cannot be AF.
\end{proof}

We say that a vertex $w$ \emph{connects to a loop} $\alpha :=
\alpha_1 \ldots \alpha_n$ if there exists a path $\gamma \in
\G^*$ with $s(\gamma)= w$ and $s(\alpha_i) \in r(\gamma)$
for some $1 \leq i \leq n$.  Note that if $w$ is a sink on a
loop (i.e. $w \in r(\alpha_i)$ for some $i$), then $w$ does
not connect to a loop.

\begin{lemma}  Let $\G$ be an ultragraph with no
singular vertices and let $A$ be the edge matrix of $\G$. 
If every vertex in $\G$ connects to a loop, then every
vertex in $\Gr(A)$ connects to a loop.
\label{con-loop}
\end{lemma}

\begin{proof}  Let $a$ be a vertex in $\Gr(A)$.  Then $a \in
\Gr(A)^0 = \G^1$.  Choose any vertex $w \in G^0$ with $w
\in r(a)$.  By hypothesis, $w$ connects to a loop $\alpha
= \alpha_1 \ldots \alpha_n$ in $\G$.  Without loss of
generality we may assume that there exists a path $\gamma =
\gamma_1 \ldots \gamma_m$ in $\G$ with $s(\gamma) = w$ and
$s(\alpha_1) \in r(\gamma)$.  Now since $A(\alpha_i,
\alpha_{i+1}) = 1$ for $1 \leq i \leq n-1$ and
$A(\alpha_n,\alpha_1)=1$, we see that there exists a loop in
$\Gr(A)$ with vertices $\alpha_1, \ldots
\alpha_n$.  Furthermore, since $A(a,\gamma_1)=1$,
$A(\gamma_i,\gamma_{i+1})=1$ for $1\leq i \leq m-1$, and
$A(\gamma_m,\alpha_1)=1$ we see that there is a path in
$\Gr(A)$ from $a$ to this loop.
\end{proof}

In \cite{Cun} Cuntz introduced the algebras $\mathcal{O}_n$
and proved that they were simple and had a property which he
called ``purely infinite".  Since that time the property of
being purely infinite has been reformulated in a number of
ways for simple $C^*$-algebras, and this has caused some
problems in deciding how to extend the notion to the
non-simple case.  In fact, various authors have used
different definitions of purely infinite for non-simple
$C^*$-algebras, and although these definitions agree in the
simple case, they are not equivalent in general.  In this
paper we shall use the definition that was used in \cite{KPR},
\cite{EL}, and \cite{BPRS}:

\begin{definition}
A $C^*$-algebra $A$ is \emph{purely infinite} if every
nonzero hereditary subalgebra of $A$ contains an infinite
projection.
\end{definition}

A competing definition is due to Kirchberg and R\o rdam
\cite{KR}: Every nonzero hereditary subalgebra of
\emph{every quotient} of $A$ contains an infinite
projection.  Note that the definition that we use is weaker
than this, but that both definitions agree in the simple
case.

\begin{theorem}  Let $\G$ be an ultragraph.  Then $C^*(\G)$
is purely infinite if and only if every loop in $\G$ has an
exit and every vertex in $\G$ connects to a loop.
\label{purely-inf}
\end{theorem}

\begin{proof}  If $\G$ contains a loop without an exit, then 
Lemma~\ref{simple-exits-nosinks} tells us that $C^*(\G)$
contains an ideal Morita equivalent to a commutative
$C^*$-algebra.  Since ideals are hereditary subalgebras this
implies that $C^*(\G)$ is not purely infinite.  

Now suppose that every loop in $\G$ contains an exit, but
that there is a vertex $v \in G^0$ that does not connect to a
loop.  Let $F^0 := \{ w \in G^0 : v \geq w \}$ and $\F^1 :=
\{ e \in \G^1 : s(e) \in F^0 \}$.  Note that $e \in \F^1$
implies $r(e) \subseteq F^0$, and thus $r$ and $s$ restrict
in such a way that we may form the ultragraph $\F := (
F^0, \F^1, r, s)$.  Let $\{s_e,p_A \}$ be the generating
Cuntz-Krieger $\G$-family.  Then Lemma~\ref{description}
implies that $\F^0 \subseteq \G^0$.  This combined with the
fact that $s(e) \in F^0$ implies $e \in \F^0$ shows that $\{
s_e, p_A : e \in \F^1, A \in \F^0 \} \subseteq C^*(\G)$ is a
Cuntz-Krieger $\F$-family.  Since $v$ does not connect to a
loop, we see that $\F$ has no loops and hence $\F$
satisfies Condition~(L).  Thus by the Cuntz-Krieger
uniqueness theorem \cite[Theorem~6.7]{Tom3}, we see
that $C^*(\F)$ is isomorphic to the subalgebra
$B := \overline{\text{span}} \{ s_\alpha p_A s_\beta^* :
\alpha, \beta \in \F^0, A \in \F^0 \}$.  We shall show that
this subalgebra is hereditary.  Let $\alpha, \beta, \gamma,
\delta \in \F^*$ and $A,B \in \F^0$.  Then for any $\mu,
\nu \in \G^*$ and $C \in \G^0$ we see from a consideration of
cases that $s_\alpha p_A s_\beta^* (s_\mu p_C s_\nu^*)
s_\gamma p_B s_\delta^*$ will have the form $s_\epsilon p_D
s_\sigma^*$ for some $\epsilon,\sigma
\in \F^*$ and $D \in \F^0$.  Since these elements span dense
subsets in $C^*(\G)$ and $B$, we see that for all $b,b' \in
B$ and $a \in C^*(\G)$ we have $bab' \in B$.  It follows from
\cite[Theorem 3.2.2]{Mur} that $B$ is hereditary.  But now,
since $\F$ does not contain any loops, Theorem~\ref{AF}
implies that $C^*(\F) \cong B$ is AF.  Hence $C^*(\G)$
cannot be purely infinite.

Conversely, suppose that $\G$ is an ultragraph in which
every loop has an exit and every vertex connects to a loop. 
Let $\F$ be a desingularization of $\G$
\cite[Definition~6.3]{Tom3}.  Then
$\F$ satisfies Condition~(L) if and only if $\G$ does, and
also every vertex in $\F$ connects to a loop if and only if
every vertex in $\G$ connects to a loop.  Since
$C^*(\G)$ is isomorphic to a full corner of $C^*(\F)$ and
because pure infiniteness is preserved by passing to corners,
it therefore suffices to prove the converse for ultragraphs with no singular vertices.

Let us therefore assume that $\G$ has no
singular vertices.  If $A$ is the edge matrix of $\G$, then
it follows from \cite[Theorem~4.5]{Tom3} that $\OA \cong
C^*(\G)$.  Now since
$\G$ satisfies Condition~(L), it follows from
\cite[Lemma~5.8]{Tom3} that
$\Gr(A)$ satisfies Condition~(L).  Also, since every vertex
in $\G$ connects to a loop, it follows from
Lemma~\ref{con-loop} that every vertex in $\Gr(A)$ connects
to a loop.  Therefore, \cite[Theorem 16.2]{EL} implies that
$\OA \cong C^*(\G)$ is purely infinite.
\end{proof}

\begin{proposition}[The Dichotomy]  Let $\G$ be an ultragraph for which $C^*(\G)$ simple.  Then
\begin{enumerate}
\item  $C^*(\G)$ is AF if $\G$ has no loops.
\item $C^*(\G)$ is purely infinite if $\G$ contains a loop.
\end{enumerate}
\end{proposition}

\begin{proof}  Since $C^*(\G)$ is simple, it follows from
Theorem~\ref{DT-simplicity} that $\G$ is cofinal and
satisfies Condition~(L).  If $\G$ has no loops, then
$C^*(\G)$ is AF by Theorem~\ref{AF}.  If $\G$ has a loop,
then every vertex connects to that loop due to cofinality,
and $C^*(\G)$ is purely infinite by Theorem~\ref{purely-inf}.
\end{proof}

%%%%%%%%%%%%%%%%%%%%%%%%%%%%%%%%%%%%%%%%%%%%%%%%%%%%%%%%%%%%%%%
\section{An ultragraph algebra that is neither an
Exel-Laca algebra nor a graph algebra}
\label{neither}
%%%%%%%%%%%%%%%%%%%%%%%%%%%%%%%%%%%%%%%%%%%%%%%%%%%%%%%%%%%%%%%

It was shown in \cite[Proposition~3.1]{Tom3} that graph
algebras are ultragraph algebras and in
\cite[Theorem~4.5]{Tom3} that Exel-Laca algebras are ultragraph algebras.  Here we show that this containment is
strict.  We provide an example of an ultragraph algebra
that is neither an Exel-Laca algebra nor an ultragraph
algebra.  

Let $A$ be the countably infinite matrix
$$
A = \left(\begin{smallmatrix}
1 & 0 & 0 & 1 & 1 & 1 & 1 &\\
0 & 1 & 0 & 1 & 1 & 1 & 1 &\\
0 & 0 & 1 & 1 & 1 & 1 & 1 &\\
1 & 0 & 0 & 1 & 0 & 0 & 0 & \cdots \\ 
0 & 1 & 0 & 0 & 1 & 0 & 0 & \\ 
0 & 0 & 1 & 0 & 0 & 1 & 0 & \\ 
0 & 0 & 0 & 1 & 0 & 0 & 1 & \\
  &   &   & \vdots & & & & \ddots \\ 
\end{smallmatrix}\right).
$$

\begin{lemma} If $A$ is as above, then $K_0(\OA) \cong 0$
and $K_1(\OA) \cong \Z \oplus \Z$.
\label{OA-K-th}
\end{lemma}

\begin{proof}  Let $I$ be the index set of $A$, and let
$\Ri$ denote the subring of $\ell^\infty (I)$ generated by
the rows $\rho_i$ of $A$ and the point masses $\delta_i$. 
If we let $A^t - I : \bigoplus_{I} \Z \rightarrow \Ri$, then
\cite[Theorem 4.5]{EL2} implies that $K_0(\OA) \cong \coker
(A^t-I)$ and $K_1(\OA) \cong \ker (A^t-I)$.  Now

$$
A^t-I = \left(\begin{smallmatrix}
0 & 0 & 0 & 1 & 0 & 0 & 0 &\\
0 & 0 & 0 & 0 & 1 & 0 & 0 &\\
0 & 0 & 0 & 0 & 0 & 1 & 0 &\\
1 & 1 & 1 & 0 & 0 & 0 & 1 & \cdots \\ 
1 & 1 & 1 & 0 & 0 & 0 & 0 & \\ 
1 & 1 & 1 & 0 & 0 & 0 & 0 & \\ 
1 & 1 & 1 & 0 & 0 & 0 & 0 & \\
  &   &   & \vdots & & & & \ddots \\ 
\end{smallmatrix}\right).
$$
Let us examine $\ker (A^t-I)$.  When $(A^t-I) (x_1,x_2,
\ldots) = \vec{0}$, then 
\begin{align*} 
x_4 &= 0 \\
x_5 &= 0 \\
x_6 &= 0 \\
x_1 + x_2 + x_3 + x_7 &= 0 \\
x_1 + x_2 + x_3 + x_8 &= 0 \\
x_1 + x_2 + x_3 + x_9 &= 0 \\
& \vdots 
\end{align*}
If $(x_1,x_2,\ldots) \in \bigoplus_I \Z$, then $x_n$ is
eventually zero, and the above equations reduce to $x_1
+x_2+x_3 = 0$ and $x_i = 0 $ for $i \geq 4$.  Hence $\ker
(A^t-I)$ is generated by $(-1,1,0,0,0,\ldots)$ and
$(-1,0,1,0,0,\ldots)$ and $\ker(A^t-I)$ has rank 2.  Thus
$K_1(\OA) \cong \ker (A^t-I) \cong \Z \oplus \Z$.

Next we shall show that $A^t-I$ maps onto $\Ri$.  Since $\Ri$
is a ring generated by $\{\rho_i,\delta_i\}$ we see that
$\Ri$ equals the collection of all sums of products of
the $\rho_i$'s and $\delta_i$'s.  But for the
matrix $A$ above, any product of the $\rho_i$'s
and $\delta_i$'s may be written as a sum of
$(0,0,0,1,1,1,\ldots)$ and the $\delta_i$'s.  Hence $\Ri =
\textrm{span}_\Z \{ (0,0,0,1,1,1,\ldots), \delta_i : i \in I
\}$.  But, $(A^t-I) \delta_{i+3} =\delta_i$ and $(A^t-I)
\delta_1 = (0,0,0,1,1,1,\ldots)$, so $A^t-I$ maps onto
$\Ri$.  Hence $K_0(\OA) \cong \coker(A^t-I) \cong 0$.
\end{proof}

For the matrix $A$ above, let $\G := (G^0, \G^1, r, s
)$ be the ultragraph $\G_A$ of Definition~\ref{edgeLG}. 
We define an ultragraph $\F$ by adding a single vertex $\{
w \}$ to $\G$ and a countable number of edges with source $w$
and range $G^0$.  More precisely, we define $\F := (F^0,
\F^1, r, s)$ by
$$ F^0 := \{w \} \cup G^0 \quad \quad \quad \quad \F^1 :=
\{e_i\}_{i=1}^\infty \cup \G^1 $$ and we extend $r$ and
$s$ to $\F^1$ by defining $s(e_i) = \{w \}$ and $r(e_i) =
G^0$ for all $1 \leq i < \infty$.

Note that $\G$ is unital because $G^0 \in \G^0$
\cite[Lemma~3.2]{Tom3}.  Since
$r(e_i) = G^0 \in \G^A$ for all $i$ we see from
Lemma~\ref{description} that $\F^0 = \{ A \cup \{ w \} : A
\in \G^0 \} \cup \G^0$.  It follows that $\F$ is also
unital.  Also note that $\G$ is transitive in the sense that
$x \geq y$ for all $x,y \in G^0$.

\begin{lemma}
Let $\F$ be the ultragraph described above and let $\Hi
:= \G^0$.  Then $\Hi$ is a saturated hereditary
subcollection of $\F^0$, the ideal $I_\Hi \triangleleft
C^*(\F)$ is Morita equivalent to $\OA$, and $C^*(\F) / I_\Hi
\cong \C$.
\label{F-facts}
\end{lemma}

\begin{proof}
Let $\{s_e,p_A\}$ be the generating Cuntz-Krieger $\F$-family
in $C^*(\F)$.  We shall first show that $I_\Hi$ is Morita
equivalent to $C^*(\G)$.  Note that since $\F^0 = \{ A \cup
\{ w \} : A \in \G^0 \} \cup \G^0$, $\{s_e,p_A\}$ restricts
to a Cuntz-Krieger $\G$-family.  Now $I_{\Hi} =
\overline{\mathrm{span}} \{ s_\alpha p_A s_\beta^* :
\alpha,\beta \in \F^* \text{ and } A \in \G^0  \}$ by
Lemma~\ref{IH}.  If we let $p: = p_{G^0}$, then $p
\in I_\Hi$ and $p_{G^0} I_\Hi p_{G^0}$ is generated by $\{s_e
,p_A : e \in \G^1 \text{ and } A \in \G^0 \}$.  Since $\G$ is
a transitive ultragraph that is not a single loop, we
see that $\G$ satisfies Condition~(L).  It then follows from
the Cuntz-Krieger Uniqueness Theorem
\cite[Theorem~6.7]{Tom3} that $C^*(\G) \cong I_\Hi$. 
All that remains to show is that $p I_\Hi p$ is a full
corner of
$I_\Hi$.  Suppose that
$J$ is an ideal in $I_\Hi$ containing $p I_\Hi p$.  Since
$p_{G^0} p_A p_{G^0} = p_A$ for all $A \in \G^0$ we see
that $\{p_A : A \in \G^0 \} \subseteq J$.  But then $J$
contains the generators of $I_\Hi$ and $J=I$.  Hence the
corner is full and $I_\Hi$ is Morita equivalent to $C^*(\G)$.

We shall now show that $C^*(\F) / I_\Hi \cong \C$.  To do
this we shall first show that $p_w \notin I_\Hi$. If it
was the case that $p_w \in I_\Hi = \overline{\text{span}}
\{ s_\alpha p_A s_\beta^* : \alpha,\beta \in \F^1, A \in \F^0
\}$, then we could find a linear combination such that 
$$\| p_w - \sum_{k=1}^n \lambda_k s_{\alpha_k}
p_{A_k} s_{\beta_k}^*  \| < 1.$$  Also since 
$$\| p_w \Big( p_w - \sum_{k=1}^n \lambda_k
s_{\alpha_k} p_{A_k} s_{\beta_k}^* \Big) \| \leq \| p_w
- \sum_{k=1}^n \lambda_k s_{\alpha_k} p_{A_k} s_{\beta_k}^* 
\|$$ we may assume that $| \alpha | \geq 1$ and $s(\alpha_k)
\in r(e)$.  Let $F$ be the (necessarily finite) set of
edges that are the initial edge of an $\alpha_i$.  Because
$w$ is an infinite emitter, it follows that $q := p_w -
\sum_{e \in F} s_es_e^*$ is a nonzero projection.  Hence 
$$\| p_w - \sum_{k=1}^n \lambda_k
s_{\alpha_k} p_{A_k} s_{\beta_k}^*  \| \geq \|
q \Big( p_w - \sum_{k=1}^n \lambda_k
s_{\alpha_k} p_{A_k} s_{\beta_k}^* \Big) \| = \| q \| =1$$
which is a contradiction.  Therefore $p_w \notin \Hi$, and 
$C^*(\F)/I_\Hi$ is generated by the projection
$p_w+I_\Hi$.  Consequently, $C^*(\F) / I_\Hi \cong \C$.
\end{proof}

\noindent The ideas in the proof of the following proposition
were suggested by Wojciech Szyma\'nski.

\begin{proposition}  The ultragraph algebra $C^*(\F)$ is
not an Exel-Laca algebra.
\end{proposition}

\begin{proof}  Recall that a character for $C^*(\F)$ is a
nonzero homomorphism $\epsilon : C^*(\F) \rightarrow
\C$.  We shall show that there is a unique character on
$C^*(\F)$.  Let $\{s_e,p_A\}$ be a generating Cuntz-Krieger
$\F$-family.

Since $C^*(\F)/I_\Hi \cong \C$ by Lemma~\ref{F-facts} we
see that the projection $\pi : C^*(\F) \rightarrow
C^*(\F)/I_\Hi$ is a character.  We shall now show that this
character is unique.  Let $\epsilon : C^*(\F) \rightarrow \C$
be a character.  Set $I = \ker \epsilon$.  Then $I$ is a
nonzero ideal and $\Hi := \{A \in \F^0 : p_A \in I \}$ is a
saturated hereditary subcollection.  Since $\G$ is
transitive, we see that $\F$ satisfies Condition~(L). 
Therefore, the Cuntz-Krieger Uniqueness Theorem
\cite[Theorem~6.7]{Tom3} implies that
$\ker \epsilon$ contains one of the $p_A$'s, and $\Hi$ is
nonempty.  Because $\Hi$ is nonempty and $\G$ is transitive,
it follows that $\G^0 \subseteq \Hi$.  Now since
$\epsilon$ is nonzero, we cannot also have $\{ w \}$ in
$\Hi$.  Therefore, $\Hi = \G^0$, and this implies that $p_v
\in I$ for all $v \in G^0$ and $s_e = s_ep_{r(e)} \in I$ for
all $e \in \F^1$.  Since $C^*(\F)$ is generated by $\{ s_e :
e \in \F^1 \} \cup \{ p_v : v \in F^0 = G^0 \cup \{ w \} \}$,
and $$\epsilon(p_w) = \epsilon(s_e) = 0 \quad \quad \text{for
all $w \in G^0$ and $e \in \F^1$}$$ we see that $\epsilon$ is
completely determined by its value on $p_w$.  Because $p_w$
is a projection, $\epsilon (p_w) =1$.  Thus $\epsilon$ is
unique.

Now if $C^*(\F)$ was an Exel-Laca algebra, then $C^*(\F)$
would be generated by an Exel-Laca family $\{S_i\}$.  Let
$\gamma$ be the gauge action on this Exel-Laca algebra. 
Because there is a unique character $\epsilon$ on
$C^*(\F)$, we see that $\epsilon \circ \gamma_z = \epsilon$
for all $z \in \T$.  Also, since $\epsilon$ is nonzero,
$\epsilon (S_i) \neq 0$ for some $i$.  Thus $\epsilon (S_i)
= \epsilon(\gamma_z ( S_i)) = z \epsilon (S_i)$ for
all $z \in \T$ which is a contradiction.
\end{proof}

\begin{proposition} The ultragraph algebra $C^*(\F)$ is
not a graph algebra.
\end{proposition}

\begin{proof}  Let $\Hi := \G^0$.  Then $\Hi$ is a saturated
hereditary subcollection of $\F^0$.  The short exact
sequence $0 \rightarrow I_\Hi \rightarrow C^*(\F)
\rightarrow C^*(\F) / I_\Hi \rightarrow 0$ induces the
following cyclic six term exact sequence for $K$-theory:
$$ \xymatrix{ K_0(I_\Hi) \ar[r] & K_0(C^*(\F)) \ar[r] &
K_0(C^*(\F) / I_\Hi) \ar[d] \\ K_1(C^*(\F) / I_\Hi) \ar[u] &
K_1(C^*(\F)) \ar[l] & K_1(I_\Hi) \ar[l] }$$
It follows from Lemma~\ref{F-facts} that $C^*(\F) /
I_\Hi \cong \C$.  Thus $K_0(C^*(\F) / I_\Hi) \cong \Z$ and
$K_1(C^*(\F) / I_\Hi) \cong 0$.  Also, Lemma~\ref{F-facts}
tells us that $I_\Hi$ is Morita equivalent to $\OA$, and it
then follows from Lemma~\ref{OA-K-th} that $K_0(I_\Hi) \cong
0$ and $K_1(I_\Hi) \cong \Z \oplus \Z$.  Thus the above exact
sequence becomes $$ \xymatrix{ 0 \ar[r] & K_0(C^*(\F)) \ar[r]
& \Z \ar[r] & \Z \oplus \Z \ar[r] & K_1(C^*(\F)) \ar[r] & 0
}$$ and $\rank K_0(C^*(\F)) < \rank K_1(C^*(\F))$.

Now we see that $F^0 = G^0 \cup \{ w \} \in \F^0$ and thus
$C^*(\F)$ is unital.  Therefore, if $C^*(\F)$ were the
$C^*$-algebra of a graph, then this graph would have to have
a finite number of vertices.  It then follows from
\cite[Theorem 3.2]{RS} that there exists an exact sequence
$$ \xymatrix{ 0 \ar[r] & K_1(C^*(\F)) \ar[r] & \bigoplus_V
\Z \ar[r] & \bigoplus_V \Z \oplus \bigoplus_W \Z \ar[r] &
K_0(C^*(\F)) \ar[r] & 0 }$$
for some finite sets $V$ and $W$.  Hence $\rank K_1(C^*(\F))
\leq \rank K_0(C^*(\F))$, which is a contradiction.
\end{proof}

\begin{corollary}  If $\F$ is the ultragraph
described above, then $C^*(\F)$ is neither an Exel-Laca
algebra nor a graph algebra.
\end{corollary}

%%%%%%%%%%%%%%%%%%%%%%%%%%%%%%%%%%%%%%%%%%%%%%%%%%%%%%%%%%%%%%%
\section{Viewing Ultragraph Algebras as Cuntz-Pimsner
Algebras}
\label{Cun-Pim}
%%%%%%%%%%%%%%%%%%%%%%%%%%%%%%%%%%%%%%%%%%%%%%%%%%%%%%%%%%%%%%%

Let $X$ be a Hilbert bimodule over a $C^*$-algebra
$\mathcal{A}$, in the sense that $X$ is a right Hilbert
$\mathcal{A}$-module with a left action of $\mathcal{A}$ by
adjointable operators.  In \cite{Pim} Pimsner described how
to construct a $C^*$-algebra $\OX$ from $X$.  These
Cuntz-Pimsner algebras have been shown to include many
classes of
$C^*$-algebras and consequently have been the subject of
much attention.  Pimsner originally showed that for
appropriate choices of $X$ and $\mathcal{A}$, the
Cuntz-Pimsner algebras included the Cuntz-Krieger algebras
\cite[\S1 Example 2]{Pim} as well as crossed products by
$\Z$ \cite[\S1 Example 3]{Pim}.  Since that time it has
also been shown that the $C^*$-algebras of graphs with no
sinks \cite[Proposition 12]{FLR} and the Exel-Laca algebras
\cite[Theorem 5]{Szy2} may be realized as Cuntz-Pimsner
algebras.

In this section we show that the $C^*$-algebras of ultragraphs with no sinks may also be realized as Cuntz-Pimsner
algebras using a construction similar to that in
\cite{Szy2}.  Let $\G = (G^0,\G^1,r,s)$ be an ultragraph
with no sinks.  Define $\mathcal{A}$ to be the
$C^*$-subalgebra of
$C^*(\G)$ generated by $\{p_A : A \in \G^0 \}$.  Note that
since the $p_A$'s commute and $\{p_A : A \in \G^0 \}$ is
closed under multiplication, it follows that $\mathcal{A} =
\overline{\textrm{span}} \{p_A : A \in \G^0 \}$.  Also let
$X := \overline{\textrm{span}} \{ s_ep_A : e \in \G^1, A \in
\G^0 \}$.  Then $X$ has a natural Hilbert
$\mathcal{A}$-bimodule structure with the right action given
by right multiplication, the left action given by left
multiplication, and the $\mathcal{A}$-valued inner product
given by
$\langle x , y \rangle_\mathcal{A} := x^*y$.

We shall let $\phi : \mathcal{A} \rightarrow \Li (X)$ denote
the map given by the left action; that is,
$\phi(a)(x):=ax$.  We shall also let $\K(X)$ denote the
compact operators on $X$ and $J(X) := \phi^{-1}(\K(X))$.

\begin{theorem}
If $X$ is the Hilbert bimodule defined above, then
$\OX$ is canonically isomorphic to $C^*(\G)$.
\end{theorem}

\begin{proof}
Using the language of \cite{FMR}, let $(k_X, k_\mathcal{A})$
be a universal Toeplitz representation of $X$ in $\OX$ which
is Cuntz-Pimsner covariant (i.e.~coisometric on $J(X)$).  We
shall show that
$\{ k_X(s_e), k_\mathcal{A}(p_A)\}$ is a Cuntz-Krieger
$\G$-family in $\OX$.

Since $k_\mathcal{A}$ is a homomorphism, we trivially have
$k_\mathcal{A}(p_Ap_B) =
k_\mathcal{A}(p_A)k_\mathcal{A}(p_B)$ and
$k_\mathcal{A}(p_{A\cup B}) = k_\mathcal{A}(p_A) +
k_\mathcal{A}(p_B) - k_\mathcal{A}(p_{A\cap B})$.  Because
$(k_X, k_\mathcal{A})$ is a Toeplitz representation we have
$k_X(s_e)^*k_X(s_e) = k_\mathcal{A} (
\langle s_e, s_e\rangle_\mathcal{A} ) = k_\mathcal{A}
(s_e^*s_e) = k_\mathcal{A} (p_{r(e)})$.  Also
$k_\mathcal{A}(p_{s(e)}) k_X(s_e) = k_X(p_{s(e)} s_e) =
k_X(s_e)$ so $k_X(s_e) k_X(s_e)^* \leq
k_\mathcal{A}(p_{s(e)})$.  Finally, if $v$ is the source of
finitely many vertices, then $p_v = \sum_{s(e)=v} s_es_e^*$
and
$\phi(p_v) = \sum_{s(e)=v} \Theta_{s_e,s_e}$.  It then
follows from the fact that $(k_X, k_\mathcal{A})$ is
Cuntz-Pimsner covariant that $k_\mathcal{A}(p_v) =
k_\mathcal{A}^{(1)}(\phi(p_v)) =
k_\mathcal{A}^{(1)}(\sum_{s(e)=v}
\Theta_{s_e,s_e}) = \sum_{s(e)=v} k_X(s_e) k_X(s_e)^*$. 
Hence $\{ k_X(s_e), k_\mathcal{A}(p_A)\}$ is a Cuntz-Krieger
$\G$-family and the universal property of
$C^*(\G)$ gives a homomorphism $\Phi : C^*(\G) \rightarrow
\OX$ with
$\Phi(s_e)=k_X(s_e)$ and $\Phi(p_A)=k_\mathcal{A}(p_A)$.  

Let $\psi : X \hookrightarrow C^*(\G)$ and $\pi
: \mathcal{A} \hookrightarrow C^*(\G)$ be the inclusion
maps.  Then
$(\psi,\pi)$ is a Toeplitz representation.  To see that
$(\psi,
\pi)$ is also Cuntz-Pimsner covariant, let $a \in
\mathcal{A}$ with
$\phi(a) \in
\K(X)$.  Then $\phi(a) = \lim \sum \lambda_k
\Theta_{x_k,y_k}$ and hence $a = \lim \sum \lambda_k
x_k^*y_k$.  But then 
\begin{align*} 
\pi^{(1)}(\phi(a)) &= \lim \sum
\lambda_k \pi^{(1)} (\Theta_{x_k,y_k}) = \lim \sum
\lambda_k \psi(x_k)^* \psi(y_k) \\
&= \lim \sum \lambda_k
x_k^*y_k = a = \pi (a).
\end{align*}
Since $(\psi,\pi)$ is a
Toeplitz representation which is Cuntz-Pimsner
covariant, the universal property of $\OX$
\cite[Proposition 1.3]{FMR} implies that there is a
homomorphism $\Phi' :
\OX \rightarrow C^*(G)$ which commutes with $(\psi,\pi)$. 
But then $\Phi'(k_\mathcal{A}(p_A)) = \pi(p_A) = p_A$ and
$\Phi'(k_X(s_e)) = \psi(s_e) = s_e$ so $\Phi$ and $\Phi'$
are inverses for each other. 
\end{proof}

\noindent In the remainder of this section we shall give a
description of $J(X) := \phi^{-1}(\K(X))$ in terms of the
ultragraph.

\begin{lemma}
If $X$ is the Hilbert bimodule of an ultragraph $\G$, then
$\phi (a) \in \K(X)$ implies $a \in \overline{\mathrm{span}}
\{ s_ep_As_f^* : e,f \in \G^1, A \in \G^0 \}$.
\label{compacts}
\end{lemma}

\begin{proof}
Since $X := \overline{\textrm{span}} \{ s_ep_A : e \in
\G^1, A \in \G^0 \}$ we have $\K(X) =
\overline{\textrm{span}} \{ \Theta_{s_ep_A,s_fp_B} : e,f \in
\G^1 A,B \in \G^0 \}$.  If $\phi(a) \in \K(X)$, then for any
$\epsilon > 0$ there exists a finite linear combination with
$\| \phi(a) - \sum \Theta_{s_ep_A,s_fp_B} \| < \epsilon$. 
Thus $$\| \phi(a)(x) - \sum \Theta_{s_ep_A,s_fp_B}(x) \| = 
\| ax - \sum s_ep_{A \cap B}s_f^* x \| <
\epsilon$$ for all $x \in X$.  Hence  $\| a - \sum s_ep_{A
\cap B}s_f^*  \| < \epsilon$ and the claim is proven.
\end{proof}

\noindent Throughout the following whenever $B,C \in \G^0$
we shall let $Q(B,C) := p_B - p_Bp_C$.

\begin{lemma}
If $A_1, \ldots, A_n \in \G^0$, then 
$$\sum_{I
\subseteq \{1,\ldots n\} }
Q(\cap_{i \in I} A_i , \cup_{i \notin I} A_i) = 1.$$
\label{sum-id}
\end{lemma}

\begin{proof}
Induct on $n$.  Multiply the formula for $n=k$ by
$p_{A_{k+1}} + (1-p_{A_{k+1}})$.
\end{proof}

\begin{lemma}
\label{mut-ortho}
If $\sum_{k=1}^n \lambda_k p_{A_k}$ is a finite linear
combination with $A_k \in \G^0$ for all $k$, then 
$$\sum_{k=1}^n \lambda_k p_{A_k} = \sum_{I \subseteq
\{1,\ldots, n \}} a_I Q(\cap_{i \in I} A_i, \cup_{i
\notin I} A_i)$$ where $a_I := \sum_{i \in I}
\lambda_i$. Consequently $\sum_{k=1}^n \lambda_k p_{A_k}$ can
be rewritten as a linear combination of mutually orthogonal
projections of the form $Q(B,C)$ with $B,C \in \G^0$.
\end{lemma}

\begin{proof}  For convenience of notation, let $N_n := \{
1, \ldots n \}$.  We shall prove the claim by induction on
$n$.  For $n=1$ the equality holds easily.  Therefore,
assume the equality is true for $n$ and we shall prove it
for $n+1$.
\begin{align*}
& \sum_{I \subseteq
N_{n+1}} a_I Q(\bigcap_{i \in I} A_i, \bigcup_{i \in N_{n+1}
\backslash I} A_i) \\
= &  \sum_{I \subseteq
N_n} a_I Q(\bigcap_{i \in I} A_i, \bigcup_{i \in N_{n+1}
\backslash I} A_i) +  \sum_{I \subseteq
N_n} (a_I+\lambda_{n+1}) Q(A_{n+1} \cap \bigcap_{i \in I}
A_i,
\bigcup_{i
\in N_n \backslash I} A_i) \\
= &  \sum_{I \subseteq
N_n} a_I p_{\cap_{i \in I} A_i} - a_I p_{\cap_{i \in I} A_i}
p_{ \cup_{i \in N_{n+1}
\backslash I} A_i} \\
& \qquad \qquad +  \sum_{I \subseteq
N_n} (a_I+\lambda_{n+1})
(p_{\cap_{i \in I} A_i} p_{A_{n+1}} -  p_{\cap_{i \in I} A_i}
p_{A_{n+1}} p_{\cup_{i
\in N_n \backslash I} A_i}) \\
= &  \sum_{I \subseteq
N_n} a_I p_{\cap_{i \in I} A_i} - a_I p_{\cap_{i \in I} A_i}
( p_{\cup_{i \in N_n \backslash I} A_i} + p_{A_{n+1}} -
p_{A_{n+1}} p_{\cup_{i \in N_n \backslash I} A_i}) \\
& \qquad \qquad + \sum_{I \subseteq
N_n} (a_I+\lambda_{n+1})
(p_{\cap_{i \in I} A_i} p_{A_{n+1}} -  p_{\cap_{i \in I} A_i}
p_{A_{n+1}} p_{\cup_{i \in N_n \backslash I} A_i}) \\
= &  \sum_{I \subseteq
N_n} a_I p_{\cap_{i \in I} A_i} - a_I p_{\cap_{i \in I} A_i}
 p_{\cup_{i \in N_n \backslash I} A_i} \\
& \qquad \qquad + \sum_{I \subseteq
N_n} \lambda_{n+1}
(p_{\cap_{i \in I} A_i} p_{A_{n+1}} -  p_{\cap_{i \in I} A_i}
p_{A_{n+1}} p_{\cup_{i \in N_n \backslash I} A_i}) \\
= &  \sum_{k=1}^n \lambda_k p_{A_k} +  \lambda_{n+1}
p_{A_{n+1}} \sum_{I
\subseteq N_n}
(p_{\cap_{i \in I} A_i}  -  p_{\cap_{i \in I} A_i} p_{\cup_{i
\in N_n \backslash I} A_i}) \\
= & \sum_{k=1}^{n+1} \lambda_k p_{A_k}
\end{align*}
where this last line follows from Lemma~\ref{sum-id}.

The final claim follows from the fact that the
terms $Q(\cap_{i
\in I} A_i, \cup_{i
\notin I} A_i)$ and $Q(\cap_{i \in J} A_i, \cup_{i
\notin J} A_i)$ are orthogonal when $I \neq J$.
\end{proof}

\begin{proposition}
If $\G$ is an ultragraph with no sinks and $X$ is the
Hilbert bimodule defined above, then $$\phi^{-1}(\K (X)) =
\overline{\mathrm{span}} \{ p_v : v \in G^0 \text{ and $v$
is not an infinite emitter} \}.$$
\end{proposition}

\begin{proof}
Let $I$ denote the right hand side of the above equation. 
If $v \in G^0$ is not an infinite emitter, then $p_v =
\sum_{s(e)=v} s_es_e^*$ and $\phi(p_v) = \sum_{s(e)=v}
\Theta_{s_e,s_e} \in \K(X)$.  Hence $I \subseteq
\phi^{-1}(\K (X))$.

To see the reverse inclusion let $a \in \mathcal{A}$ and
$\phi(a) \in \K (X)$.  Choose $\epsilon > 0$.  Since
$\mathcal{A} =
\overline{\textrm{span}} \{p_A : A \in \G^0
\}$, Lemma~\ref{mut-ortho} implies that there exists a
finite linear combination $\sum_{k=1}^n \lambda_k
Q(B_k,C_k)$ with the $Q(B_k,C_k)$'s mutually orthogonal and
$\| a - \sum_{k=1}^n \lambda_k Q(B_k,C_k) \| < \epsilon /
2$.  Define $s^{-1}(B \backslash C) := \{ e \in \G^1 : s(e)
\in B \backslash C \}$, and let $S_1 := \{ k : |s^{-1}(B_k
\backslash C_k)| < \infty \}$ and $S_2 := \{ k : |s^{-1}(B_k
\backslash C_k)| = \infty \}$.  Let $\lambda_{k'}:= \max \{
|\lambda_k| : k \in S_2 \}$.  Because $\phi(a) \in \K (X)$,
we know from Lemma~\ref{compacts} that $a \in
\overline{\mathrm{span}} \{ s_ep_As_f^* : e,f \in \G^1, A \in
\G^0 \}$.  Thus for every $\epsilon' > 0$ we may find a
finite linear combination $\sum_{j=1}^m \mu_j s_{e_j} p_{A_j}
s_{f_j}^*$ such that $\| a - \sum_{j=1}^m \mu_j s_{e_j}
p_{A_j} s_{f_j}^* \| < \epsilon'$.

Now since $|s^{-1}(B_{k'} \backslash C_{k'}) | = \infty$,
there exists an edge $g$ such that $s(g) \in B_{k'}
\backslash C_{k'}$ and $g$ is not equal to any of the
$f_j$'s.  Since the
$Q(B_k,C_k)$'s are mutually orthogonal, we have that
$s(g) \notin B_{k} \backslash C_{k}$ for all $k \neq k'$. 
Hence

\begin{align*} | \lambda_{k'} | & = \|
\lambda_{k'} Q(B_{k'},C_{k'}) \| \\ 
& = \| \sum_{k=1}^n
\lambda_k Q(B_k,C_k) s_g -  \sum_{j=1}^m \mu_j s_{e_j}
p_{A_j} s_{f_j}^* s_g \| \\
& \leq \| \sum_{k=1}^n \lambda_k Q(B_k,C_k) -  \sum_{j=1}^m
\mu_j s_{e_j} p_{A_j} s_{f_j}^*\| \ \|s_g \| \\
& \leq \| \sum_{k=1}^n \lambda_k Q(B_k,C_k) -  a \| + \|
a - \sum_{j=1}^m \mu_j s_{e_j} p_{A_j} s_{f_j}^*\| \\
& \leq \| \sum_{k=1}^n \lambda_k Q(B_k,C_k) -  a \| +
\epsilon'.
\end{align*}
Since this inequality holds for all $\epsilon' > 0$ we have
$ | \lambda_{k'} | \leq   \| \sum_{k=1}^n \lambda_k
Q(B_k,C_k) -  a \| < \epsilon / 2$.  Thus 
\begin{align*} \| a - \sum_{k \in S_1} \lambda_k Q(B_k,C_k)
\| & \leq \| a - \sum_{k=1}^n \lambda_k Q(B_k,C_k) \| + \|
\sum_{k \in S_2}
\lambda_k Q(B_k,C_k) \| \\
& < \frac{\epsilon}{2} + | \lambda_{k'} | \leq
\frac{\epsilon}{2} + \frac{\epsilon}{2} = \epsilon.
\end{align*}

Now for every $k \in S_1$ we have that $s^{-1}(B_k
\backslash C_k) < \infty$.  Since $\G$ has no sinks, this
implies that $B_k \backslash C_k$ is the union of a finite
number of vertices that emit finitely many edges.  Since
$B_k \backslash C_k$ is finite \cite[Lemma 4.2]{Tom3}
implies that $p_{B_k \backslash C_k} = \sum_{v \in B_k
\backslash C_k} p_v$.  Furthermore, since $B_k \backslash
C_k$ is finite, it is an element of $\G^0$ and the equality
$p_{B_k} = p_{B_k \backslash C_k} + p_{B_k \cap C_k} -
p_\emptyset$ shows that $Q(B_k,C_k) = p_{B_k \backslash
C_k} = \sum_{v \in B_k \backslash C_k} p_v$.  Since the
$p_v$'s all emit finitely many edges and since $\epsilon$
was arbitrary, the above shows that $a \in
\overline{\mathrm{span}} \{ p_v : v \in G^0 \text{ and $v$
is not an infinite emitter} \}$ and hence $\phi^{-1}(\K (X))
\subseteq I$.
\end{proof}

\begin{corollary}
If $\G$ is an ultragraph with no sinks, then
$\phi^{-1}(\K(X)) \cong C_0(T)$ where $T := \{ v \in G^0 :
v \text{ emits finitely many edges} \}$ has the discrete
topology.
\end{corollary}

\begin{proof}
Let $\delta_v \in C_0(T)$ denote the point mass at $v$. 
Then the map $\delta_v \mapsto p_v$ extends to an
isomorphism from $C_0(T)$ onto $\overline{\text{span}}
\{p_v : v \text{ emits finitely many edges} \}$. 
\end{proof}

\begin{remark}
If $G$ is a graph, then the $C^*$-algebra $\mathcal{A}$ in
the graph bimodule \cite[Example 1.2]{FR} is equal to
$C_0(G^0)$ for the discrete space $G^0$.  For an ultragraph
$\G$ the $C^*$-algebra
$\mathcal{A}$ arising in the ultragraph bimodule is the
$C^*$-algebra generated by $\{p_A : A
\in \G^0 \}$.  Since the $p_A$'s commute $\mathcal{A}$ is
commutative and $\mathcal{A} \cong C_0(X)$ for some locally
compact space $X$.  Furthermore, since $\mathcal{A}$ is
generated by projections the space $X$ must be totally
disconnected.  However, in general $X$ need not be
discrete.  Despite this, the above corollary shows that the
ideal $\phi^{-1}(\K(X))$ corresponds to $C_0(T)$ for some
discrete open set $T
\subseteq X$. 
\end{remark}

\end{document}